\documentclass[reqno]{amsart}  

\usepackage{amsmath,amsfonts,amssymb,mathtools}
\usepackage{amsthm}
\usepackage[hidelinks]{hyperref}
\usepackage{xcolor}
\usepackage{graphicx}
\usepackage{pgfplots}
\usepackage{tikz}
\usepgfplotslibrary{groupplots}
\pgfplotsset{compat=1.14}
\newcommand{\legendline}[1]{%
	\protect(\protect\tikz[baseline=-0.5ex]
	{\protect\draw[#1] (0,0) -- (6mm,0);}\protect)%
}
\usepackage{standalone}
\usepackage{booktabs,tabularx}
\usepackage{siunitx}
\newcommand{\R}{\mathbb R}
\newcommand{\C}{\mathbb C}
\newcommand{\N}{\mathbb N}
\newcommand{\B}{\mathcal B}

\newcommand{\columnsselector}{\mathfrak g}
\newcommand{\intv}{\mathcal I}

\DeclareMathOperator{\col}{\operatorname{col}}
\DeclareMathOperator{\rank}{\operatorname{rank}}
\newcommand{\tddt}{\tfrac{\mathrm d}{\mathrm dt}}
\DeclareMathOperator{\lag}{\mathfrak l}

\newtheorem{theorem}{Theorem}
\newtheorem{lemma}[theorem]{Lemma}
\newtheorem{proposition}[theorem]{Proposition}
\newtheorem{definition}[theorem]{Definition}
\newtheorem{corollary}[theorem]{Corollary}
\newtheorem{remark}[theorem]{Remark}
\newtheorem{assumption}{Assumption}

\newcommand{\timevar}{t}                           
\newcommand{\timewindow}{T_{\text{diff}}}                        

\newcommand{\diffkernelsymbol}{g}                     
\newcommand{\diffkernel}[1]{\diffkernelsymbol^{(#1)}}               

\newcommand{\jacoparam}{\alpha_{\text{diff}}}                    
\newcommand{\jacoparambeta}{\beta_{\text{diff}}}                 
\newcommand{\polyorder}{N_{\text{diff}}}                         
\newcommand{\jacopoly}[1]{P_{#1}^{(\jacoparam,\jacoparambeta)}} 
\newcommand{\derivorder}{n}                        

\title{
A data-based image representation for continuous-time LTI systems
}

\author[A.~Othmane]{Amine Othmane${}^\dagger$}
\address{Systems Modeling and Simulation, Saarland University, Saarbrücken, Germany}
\email{\href{mailto:amine.othmane@uni-saarland.de}{amine.othmane@uni-saarland.de}}
\author[P.~Schmitz]{Philipp Schmitz${}^\dagger$}
\address{Optimization-based Control Group, Technische Universität Ilmenau, Ilmenau, Germany}
\email{\href{mailto:philipp.schmitz@tu-ilmenau.de}{philipp.schmitz@tu-ilmenau.de}}
\author[K.~Worthmann]{Karl Worthmann}
\address{Optimization-based Control Group, Technische Universität Ilmenau, Ilmenau, Germany}
\email{\href{mailto:karl.worthmann@tu-ilmenau.de}{karl.worthmann@tu-ilmenau.de}}
\author[K.~Flaßkamp]{Kathrin Flaßkamp}
\address{Systems Modeling and Simulation, Saarland University, Saarbrücken, Germany}
\email{\href{mailto:kathrin.flasskamp@uni-saarland.de}{kathrin.flasskamp@uni-saarland.de}}

\thanks{${}^\dagger$Both authors contributed equally to this work. A.~Othmane is grateful for the support from the German Federal Ministry of Research, Technology and Space (No. 05M22TSA). P.~Schmitz is grateful for the support from the Carl Zeiss Foundation (VerneDCt -- No. 2021-10-003).}

\begin{document}

\begin{abstract}
We derive a numerically stable method to compute an image representation of an unknown linear system only from data, leveraging a continuous-time version of Willems et al.'s fundamental lemma. To this end, we use derivatives approximated by algebraic differentiators. Our novel image representation avoids solving differential-algebraic equations and significantly reduces computational complexity by eliminating redundant degrees of freedom corresponding to the number of unknown quantities to be identified. Simulation results confirm the effectiveness of the proposed approach, even in the presence of severe measurement disturbances.
\end{abstract}

\maketitle

\section{Introduction}

design of controllers and feedback laws directly from data, commonly referred to as direct data-driven control (see, e.g.,~\cite{dorfler2022bridging}), has attracted considerable research interest in recent years. A key result in this field is the fundamental lemma by Willems et al.~\cite{Willems2005}. This lemma states that all finite-length trajectories of a controllable discrete-time linear time-invariant (LTI) system can be parameterized using input--output data arranged in a Hankel matrix, provided suitable persistency of excitation conditions are met. This result has enabled numerous applications in data-driven control, e.g.\ in~\cite{Persis2020} and~\cite{Coulson2019}; see \cite{Markovsky2021} and \cite{Faulwasser2023} for comprehensive overviews.

Several continuous-time extensions have been developed. The approach in \cite{Lopez_2022} uses Hankel-like structures from sampled input--state trajectories. Input--output versions that avoid state measurements were proposed in \cite{Rapisarda2023} using orthogonal polynomials, and in \cite{Lopez2024} and \cite{Schmitz2024a, Schmitz2026} using input--output jets arranged in Hankel-like and Gramian structures, respectively. The orthogonal polynomial approach in \cite{Rapisarda2023} yields in practice an approximate fundamental lemma, whereas \cite{Lopez2024} and \cite{Schmitz2024a, Schmitz2026} provide exact characterizations but require solving a differential-algebraic equation (DAE). Moreover, these input--output methods depend on higher-order derivatives of an informative trajectory, which are typically unknown and must be estimated from noisy measurements. Computing derivatives from noisy data is challenging, as it constitutes an ill-posed problem.

In this paper, we address the two main challenges discussed above, namely the need to solve DAEs and the estimation of higher-order derivatives from noisy measurements. Based on the result in \cite{Schmitz2024a}, we first derive an image-based, data-driven formulation for controllable continuous-time LTI systems, cf.~\cite{Willems1997}. This representation circumvents the need to solve DAEs and reduces the complexity of the data-based trajectory parameterization by removing redundant degrees of freedom, so that fewer unknowns need to be identified. Second, we incorporate algebraic differentiators introduced in \cite{mboup2009a} to reliably estimate the required derivatives from measured data. These differentiators provide an efficient and robust approach, as evidenced by the diverse applications summarized in \cite{othmane2021b}, which also contains tuning guidelines. Moreover, an open-source software toolbox accompanied by a practical introduction can be found in \cite{othmane2023Tool}. We show that the estimated signals preserve the persistency of excitation property required by the fundamental lemma. We restrict attention to single-output systems and develop the method in this setting. Simulation results demonstrate the effectiveness of the approach, even in the presence of substantial measurement disturbances.

The paper is organized as follows. Section~\ref{sec:setting} formalizes the problem setting and recalls the concepts of persistency of excitation and the continuous-time fundamental lemma. In Section~\ref{sec:ddimrep}, we develop a method to construct an image representation from noisy data. Section~\ref{sec:experiments} illustrates the theoretical findings using numerical simulations before conclusions are drawn in Section~\ref{sec:conclusions}. Supplementary material on algebraic differentiators, including the definition of suitable filter kernels, is provided in Appendix~\ref{app:algdiff}.

\textbf{Notation.} Given an interval $\intv\subset \R$ and $d\in\N\setminus\{0\}$, the space of equivalence classes of square integrable functions from $\intv$ to $\R^d$ is denoted by $L^2(\intv,\R^d)$. The associated Sobolev space of order $k\in\N$ is denoted by $H^k(\intv, \R^d)$. The convolution of two integrable functions is defined as $(f* g)(t):=\int_{\R} f(s)g(t-s)\,\mathrm ds$, where functions defined on an interval $\intv\subset \R$ are extended by zero outside of $\intv$. For a sufficiently smooth function $x:\mathcal I\to\R$, we denote by $x^{(i)}$
its $i$-th derivative, with the convention ${x^{(0)} = x}$; we also write $\dot x = x^{(1)}$. For matrices or vectors $A_1,\dots, A_n$ of compatible size we write
$\col(A_1,\dots,A_n) := \begin{bmatrix}
            A_1^\top & \dots & A_n^\top
        \end{bmatrix}{}^\top$.  
The notation $\mathrm{Unif}([-a, a])$ with $a>0$, refers to the uniform distribution on the interval $[-a, a]$.

\section{Setting and problem description}
\label{sec:setting}

In this section, we introduce the problem setting and recall the continuous-time fundamental lemma.

\subsection{Linear time-invariant systems}
We focus on continuous-time LTI systems of the form
\begin{equation}
    \label{eq:sys}
    \begin{split}
        \dot x(t) &= Ax(t) + Bu(t)\\
        y(t) &= Cx(t) + Du(t)
    \end{split}
\end{equation}
on an open and finite time interval $\intv$ with system matrices $A\in\R^{n\times n}$, $B\in\R^{n\times m}$, $C\in\R^{p\times n}$ and $D\in\R^{p\times m}$. 
At time $t\in\intv$, the state, the input, and the output of the system are given by $x(t)\in\R^n$, $u(t)\in \R^m$ and $y(t)\in\R^p$, respectively. As the set of valid input--output trajectories corresponding to system~\eqref{eq:sys}, we consider the \emph{behavior}
\begin{equation}
    \B :=\left\{\col(u,y)\in L^2(\intv, \R^{m+p})\,\middle|\, \begin{aligned}& \exists\, x\in H^1(\intv, \R^n) \text{ s.t.}\\&\text{\eqref{eq:sys} holds a.e.\ in }\intv\end{aligned}\right\}.
\end{equation}
The behavior $\B$ is a closed linear subspace in $L^2(\intv, \R^{m+p})$, cf.\ \cite[Lem.~2 and Lem.~3]{Schmitz2024a}.

\begin{assumption}
\label{as:contrl_obs}
    Assume system \eqref{eq:sys} is controllable and observable.
\end{assumption}
For algebraic conditions equivalent to controllability and observability, 
see, e.g.,\ \cite[Section~5]{Willems1997}. Then
controllability is equivalent to the existence of an \emph{image representation} of the behavior, see \cite[Theorem~6.6.1]{Willems1997} and \cite[Lemma~6]{Schmitz2024a} for the setting of closed $L^2$-subspaces. Accordingly, there is a polynomial matrix ${M\in \R^{(m+p)\times m}[s]}$ with $M(s)=\sum_{k=0}^{n+1} M_k s^k$ such that
\begin{equation}
\label{eq:imrep}
    \B = \left\{w\in L^2(\intv, \R^{m+p})\,\middle|\, \begin{aligned} & \exists\, \ell\in L^2(\intv, \R^m)\text{ s.t.}\\ & w=M(\tddt) \ell \text{ a.e. in } \intv\end{aligned}\right\}.
\end{equation}
The variable $\ell$ in \eqref{eq:imrep}, called the \emph{latent variable}, is a differentially flat output of system~\eqref{eq:sys}.
The equation $M(\tddt)\ell = w$ should be interpreted in the weak sense, i.e.\
\begin{equation}
\label{eq:weak}
    \int_\intv w^\top(s) \phi(s)\mathrm ds = \sum_{k=0}^{n+1} (-1)^k\int_\intv \ell(s)^\top M_k^\top\phi^{(k)}(s)\,\mathrm ds
\end{equation}
for every infinitely differentiable function $\phi:\intv\rightarrow \R^m$ with compact support. 
For any $k$-times continuously differentiable $\ell:\intv\rightarrow\R^m$, the pointwise equation $M(\tddt)\ell = w$ is equivalent to \eqref{eq:weak}, which follows by integration by parts. 

An important quantity associated with $\B$ is the \emph{lag} $\lag(\B)$ (see~\cite{Willems1997}), which coincides with the observability index of \eqref{eq:sys}, i.e.\ the smallest integer $k\geq 1$ such that
\begin{equation*}
    \rank\begin{bmatrix}
        C^\top & (CA)^\top & \dots & (CA^{k-1})^\top
    \end{bmatrix} = n.
\end{equation*}

\subsection{A data-driven representation}
In the following, we recall the concept of \textit{persistency of excitation} and a version of the fundamental lemma for continuous-time LTI systems as developed in \cite{Schmitz2024a}.

\begin{definition}
    A function $u:\intv\rightarrow\mathbb R^{m}$ is called persistently exciting of order $k$, $k \in \mathbb{N} \setminus \{0\}$, if $u\in H^{k-1}(\intv, \mathbb R^{m})$ and the set of the $k\cdot m$ scalar components of $\col(u,\dots, u^{(k-1)})$ is linearly independent in $L^2(\intv,\mathbb R)$.
\end{definition}

\begin{lemma}[Fundamental lemma]
    \label{lem:FL}
    Let Assumption \ref{as:contrl_obs} hold and consider a trajectory  ${\col(\bar u,\bar y)\in\B}$  such that $\bar{u}$ is persistently exciting of order $L+n$ for some $L\geq \lag(\B)+1$. Define the Gramian matrix
    \begin{equation}
    \label{eq:gamma}
        \Gamma := \int_\intv \overline W(s)\overline W^\top\!(s)\,\mathrm ds\;\in\mathbb
        R^{L(m+p)\times L(m+p)}
    \end{equation}
    for $\overline W = \col(\bar u,\dots, \bar u^{(L-1)},\bar y,\dots, \bar y^{(L-1)})$.
    Then, ${\col(u,y)\in H^{L-1}(\intv, \R^{m+p})}$ is an element of $\B$ if and only if there exists $\columnsselector\in L^2(\intv,\mathbb R^{L(m+p)})$ satisfying
    \begin{equation}
        \label{eq:FL}
        \col(u,\dots, u^{(L-1)}, y,\dots, y^{(L-1)}) = \Gamma \columnsselector.
    \end{equation}
    Moreover, $\rank\Gamma=Lm+n$.
\end{lemma}
\begin{remark}
    \label{rem:replace_gamma}
    In \eqref{eq:FL}, $\Gamma$ may be replaced by any matrix $\widetilde \Gamma$ with $\operatorname{im}\widetilde\Gamma = \operatorname{im}\Gamma$, with $\mathfrak g$ taken in the corresponding $L^2$-space.
\end{remark}

The fundamental lemma allows for a complete description of the system's behavior by means of the data matrix $\Gamma$. This can be used, e.g.,\ for data-based system simulation or in developing direct data-based control algorithms.

Unlike the discrete-time fundamental lemma by~\cite{Willems2005}, where finite-horizon trajectories are obtained directly via a matrix–vector product, the continuous-time case in \eqref{eq:FL} requires solving a DAE. 
Using the decomposition
\begin{equation}
\label{eq:gamma_decomp}
    \Gamma=\col(\Gamma_{\bar{u}^{(0)}},\dots, \Gamma_{\bar{u}^{(L-1)})}, \Gamma_{\bar y^{(0)}},\dots, \Gamma_{\bar{y}^{(L-1)}})
\end{equation}
with ${\Gamma_{\bar u^{(k)}}\in \R^{m\times L(m+p)}}$, ${\Gamma_{\bar{y}^{(k)}}\in\R^{p\times L(m+p)}}$ for $k\in\{0,\dots, L-1\}$, we obtain the equivalence of \eqref{eq:FL} to $\col(u,y) = \col(\Gamma_{\bar{u}^{(0)}},\Gamma_{\bar{y}^{(0)}})\columnsselector$ together with
\begin{equation}
\label{eq:dae}
    \frac{\mathrm d}{\mathrm dt}(E_0 \columnsselector) = A_0 \columnsselector,
\end{equation}
where
\begin{equation*}
    \begin{aligned}
        E_0 &= \col(\Gamma_{\bar{u}^{(0)}},\dots, \Gamma_{\bar{u}^{(L-2)}},
        \Gamma_{\bar{y}^{(0)}},
        \Gamma_{\bar{y}^{(L-2)}}),\\
        A_0 &= \col(\Gamma_{\bar{u}^{(1)}},\dots, \Gamma_{\bar{u}^{(L-1)}},
        \Gamma_{\bar{y}^{(1)}},\dots,
        \Gamma_{\bar{u}^{(L-1)}}).
    \end{aligned}
\end{equation*}

\section{Image representation from noisy data}
\label{sec:ddimrep}
In the following, we present a procedure for deriving an image representation \eqref{eq:imrep} from the data-driven model \eqref{eq:FL}, starting from nominal data and then extending the approach to estimated derivatives from noisy measurements. Our representation eliminates the need to solve the DAE~\eqref{eq:dae} and, at the same time, removes redundancies in the data-based description. While $\ell\in L^2(\intv,\R^m)$ in \eqref{eq:imrep}, $\mathfrak g$ is required to belong to $L^2(\intv,\R^{L(m+p)})$ in \eqref{eq:FL}, introducing redundant degrees of freedom. The elimination of this redundancy constitutes one of the central contributions of this work.

\subsection{Unimodular embedding}
A key step, based on nominal data, is to embed the pencil $sE_0-A_0$ (associated with the DAE~\eqref{eq:dae}) into a unimodular polynomial matrix. A polynomial matrix $G$ is called \emph{unimodular} if it is square and invertible over the ring of polynomial matrices.

The required embedding therefore consists of finding a polynomial matrix $K$ such that $G$ with
\begin{equation}
\label{eq:unimod}
    G(s)=\begin{bmatrix} sE_0^\top - A_0^\top& K^\top(s)\end{bmatrix}^\top
\end{equation}
is unimodular.
\begin{proposition}[{{\cite[Cor.~1]{Beelen1988}}}]
\label{prop:unimod_emb}
    A matrix pencil $sE_0-A_0$ has a unimodular embedding~\eqref{eq:unimod} if and only if it has full row rank for all $s\in\C$. In this case, there exists a unimodular embedding with constant $K$.
\end{proposition}

In the following, we restrict our attention to the single-output case ($p=1$), which simplifies the analysis. In this case, the lag $\lag(\B)$ and the observability index both equal the system order $n$.

\begin{corollary}
\label{cor:p=1}
    Assume ${p=1}$ and let the assumptions of Lemma~\ref{lem:FL} hold with $L=\lag(\B)+1=n+1$. Then the matrix pencil $sE_0-A_0$ associated with~\eqref{eq:dae} has a unimodular embedding~\eqref{eq:unimod} with constant $K$.
\end{corollary}

\begin{proof}
    Define $\mathcal O_k = \col(C,\dots, CA^{k-1})\in\R^{k\times n}$,
    \begin{equation}
        \R^{k\times km}\ni\mathcal T_k = \begin{cases}
            D & k = 1,\\
            \left[\begin{smallmatrix}
                D \\
                \mathcal O_{k-1}B & \mathcal T_{k-1}
            \end{smallmatrix}\right], & k\in\{2,\dots, n+1\}
        \end{cases}
    \end{equation}
       and the full column rank matrix
    \begin{equation}
        \mathcal P = \begin{bmatrix}
            \mathcal O_{n+1} & \mathcal T_{n+1}\\
            0 & I_{(n+1)m}
        \end{bmatrix}\in\R^{(n+1)(m+1) \times (n+(n+1)m)}.
    \end{equation}
    Let $\Gamma_0 = \int_{\intv} W(t)W^\top(t)\,\mathrm dt\;\in\R^{(n+(n+1)m)\times (n+(n+1)m)}$ with $W=\col(\bar x, \bar u,\dots, \bar u^{n})$. By Proposition~21 in \cite{Schmitz2024a}, $\Gamma_0$ is positive definite. Then $sE_0-A_0$ equals
    \begin{equation}
        \left(s\begin{bmatrix}
            \mathcal O_n & \mathcal T_n & 0_{n\times m}\\
            0 & I_{nm} & 0_{nm\times m}
        \end{bmatrix}
        - \begin{bmatrix}
            \mathcal O_n A & \mathcal O_n B & \mathcal T_n\\
            0_{nm\times n} & 0_{nm\times m} & I_{nm}
        \end{bmatrix}\right)\Gamma_0 \mathcal P^\top.
    \end{equation}
    Note that $\Gamma_0 \mathcal P^\top$ has full row rank.
    With the matrices
    \begin{equation}
        S = \begin{bmatrix}
            \mathcal O_n^{-1} & -\mathcal O_n^{-1}\mathcal T_n\\
            0 & I_{nm}
        \end{bmatrix},\quad L(s) = \left[\begin{smallmatrix}
            s I_m & -I_m\\
            & \ddots & \ddots\\
            & & sI_m & -I_m\\
        \end{smallmatrix}\right],
    \end{equation}
    where $S$ is invertible and $L(s)\in \mathbb R[s]^{nm \times (n+1)m}$,
    we find
    \begin{equation}
        S(sE_0-A_0) = \begin{bmatrix}
            sI_n - A & -B & 0_{n\times nm}\\
            0 & \multicolumn{2}{c}{L(s)}
        \end{bmatrix} \Gamma_0\mathcal P^\top.
    \end{equation}
    The Hautus lemma for controllability together with the structure of $L(s)$ and the full row rank of $\Gamma_0 \mathcal P^\top$ implies that $S(sE_0-A_0)$ and, therefore, $sE_0-A_0$ has full row rank for every $s\in\C$. Proposition~\ref{prop:unimod_emb} yields the assertion.
\end{proof}

Let $w\in\B\cap H^{L-1}(\intv, \R^{m+p})$. Then by Lemma~\ref{lem:FL} there is $\columnsselector\in L^2(\intv, \R^{L(m+p)})$ such that \eqref{eq:FL} holds. Given a unimodular embedding $G$ of $sE_0 - A_0$, one has
\begin{equation}
    G(\tddt)\columnsselector = \begin{bmatrix}
        \tddt E_0 -A_0\\ K
    \end{bmatrix} \columnsselector = \begin{bmatrix}
        0\\\ell
    \end{bmatrix},
\end{equation}
where $\ell := K\columnsselector$. By inverting $G$, one arrives at an image representation~\eqref{eq:imrep},
\begin{equation}
\label{eq:imrep2}
    w = \begin{bmatrix}
        \Gamma_{\bar u}\\
        \Gamma_{\bar y}
    \end{bmatrix}G^{-1}(\tddt) \begin{bmatrix}
        0\\ \ell
    \end{bmatrix} = \begin{bmatrix}
        \Gamma_{\bar u}\\
        \Gamma_{\bar y}
    \end{bmatrix}G^{-1}(\tddt) \begin{bmatrix}
        0\\ I
    \end{bmatrix} \ell.
\end{equation}

\subsection{Staircase form}
\label{sec:staircase}
To construct a unimodular embedding, we transform the pencil $sE_0-A_0$ corresponding to \eqref{eq:dae} into \emph{staircase form} using the algorithm described in~\cite{oarua1997improved}. The staircase form reveals the Kronecker canonical structure of the pencil. 

Under the conditions of Corollary~\ref{cor:p=1}, there are orthogonal matrices $Q$, $Z$ such that
\begin{equation*}
\label{eq:staircase}
\begin{split}
&Q^\top(sE_0-A_0)Z = sE_1-A_1=\\[0.5em]
    &\begin{bmatrix}
        -A_{1,1} & sE_{1,2}-A_{1,2} & \cdots & sE_{1,k}-A_{1,k}\\[1em]
                & \ddots & \ddots & \vdots\\[1em]
                &        & -A_{k-1,k-1} & sE_{k-1, k}-A_{k-1, k}\\[1em]
                &        &         & -A_{k,k}
    \end{bmatrix},
\end{split}
\end{equation*}
where $A_{i,i}$, $i=1,\dots,k+1$, are upper-triangular with full row rank, while $E_{i,i+1}$, $i=1,\dots, k$, have full column rank, see \cite[Prop.~1]{oarua1997improved} and \cite{Beelen1988}.
A numerically stable algorithm to calculate the staircase form is described in \cite{oarua1997improved}.  

The pencil $(sE_1-A_1)$ can be easily embedded in a unimodular pencil. Corresponding to each $A_{i,i}$, find a matrix $C_i$ such that
\begin{equation*}
\tilde A_{i,i}:=\begin{bmatrix}
        A_{i,i}\\ K_i
    \end{bmatrix}
\end{equation*}
is upper-triangular and invertible, i.e.,\ let $K_i$ be the lower part of a truncated identity matrix. Consider the embedding of $(sE_1-A_1)$ into the pencil
$
    (sE_2-A_2) = \begin{bmatrix}
        sE_1^\top - A_1^\top&
        \tilde K^\top
    \end{bmatrix}{}^\top,
$ with $\tilde K = -\operatorname{diag}(K_1,\dots,K_{k})$.
Let $P$ be the permutation matrix such that $(sE_3-A_3)=P(E_2-sA_2)$, which is of the form
\begin{equation*}
sE_3-A_3=
\begin{bmatrix}
        -\tilde A_{1,1} & \left[\begin{smallmatrix}sE_{1,2}- A_{1,2}\\ 0\end{smallmatrix}\right] & \cdots & \left[\begin{smallmatrix}sE_{1,k}- A_{1,k}\\ 0\end{smallmatrix}\right]\\[1.5em]
                & \ddots & \ddots & \vdots\\[1.5em]
                &        & -\tilde A_{k-1,k-1} & \left[\begin{smallmatrix}sE_{k-1,k}- A_{k-1,k}\\ 0\end{smallmatrix}\right]\\[1.5em]
                &        &         & -\tilde A_{k,k}
    \end{bmatrix}.
\end{equation*}
Note that $A_3$ is upper-triangular and invertible. The matrix $E_3$ is nilpotent and the same holds, therefore, for $N:=E_3 A_3^{-1}$. Let $\eta\in\mathbb N$ be the nilpotency index of $N$, i.e.,\ $N^{\eta-1}\neq 0$ and $N^\eta=0$. Then the inverse of the pencil $(sE_3-A_3)$ is given by the polynomial matrix
\begin{equation}
\label{eq:invform}
    (sE_3-A_3)^{-1} = A_3^{-1} (sN-I)^{-1} = -A_3^{-1} \sum_{j=0}^{\eta-1} s^j N^j
\end{equation}
and, hence, $sE_3-A_3$ is unimodular.

\subsection{Image representation}
The next theorem combines the previous findings.
\begin{theorem}
    \label{thm:ddimrep}
    Under the assumption of Corollary~\ref{cor:p=1}, consider the decomposition~\eqref{eq:gamma_decomp}, and let $sE_3-A_3$ be the pencil obtained from the staircase algorithm together with the unimodular embedding using the orthogonal matrices $Z$, $Q$ and $P$ from Section~\ref{sec:staircase}.  Then the matrix $M(s)$ given via
    \begin{equation}
    \label{eq:ddimrep}
        M(s) = -\begin{bmatrix}
            \Gamma_{\bar{u}}\\\Gamma_{\bar{y}}
        \end{bmatrix} Z A_3^{-1}\sum_{j=0}^{n-1} s^j \bigl(E_3 A_3^{-1}\bigr)^{j} P\begin{bmatrix}
        0 \\ I
    \end{bmatrix}
    \end{equation}
    defines an image representation of $\B$ of the form~\eqref{eq:imrep}.
\end{theorem}
\begin{proof}
    A straightforward calculation shows that
\begin{equation}
    G(s):=\begin{bmatrix}
        sE_0-A_0\\
        \tilde KZ^\top
    \end{bmatrix} = \begin{bmatrix}
        Q\\ & I
    \end{bmatrix}P^\top(sE_3-A_3) Z^\top
\end{equation}
is a unimodular embedding of the pencil $sE_0-A_0$ corresponding to \eqref{eq:dae} with inverse
\begin{equation}
\label{eq:unimod_inv}
    G(s)^{-1} = - Z \left(A_3^{-1}\sum_{j=0}^{n-1} s^j \bigl(E_3 A_3^{-1}\bigr)^{j}\right) P\begin{bmatrix}
        Q^\top& \\& I
    \end{bmatrix}.
\end{equation}
 The assertion follows together with \eqref{eq:imrep2}.
\end{proof}

Involving only orthogonal transformations, the calculation of the staircase form is numerically stable, especially in the presence of noisy or perturbed data. The inversion of $A_3$ during the computation of the image representation, however, is highly sensitive to data corruption; a strategy is discussed in Section~\ref{sec:experiments}.

\subsection{Gramians from noisy data}

A major challenge for the application of the method discussed arises from the fact that the data matrix \(\Gamma\) contains higher-order derivatives of the input and output signals contained in \(\operatorname{col}(\bar{u}, \bar{y})\), which are typically corrupted by noise. However, the algebraic differentiators in~\cite{mboup2009a} offer an efficient and robust solution to this problem. 
A concise overview of algebraic differentiators  is provided in the appendix.

In the following, we show the invariance of persistency of excitation under convolution with a suitable filter kernel.

\begin{assumption}
\label{as:1}
    Let $g:\R\rightarrow\R$ be a bounded, $(L-1)$-times continuously differentiable function satisfying the commutation property  $\tfrac{\mathrm d^k}{\mathrm dt^k}(f*g) = (f^{(k)}*g) = (f*g^{(k)})$ for every $f\in H^{L-1}(\intv, \R)$ and all $k=1,\dots,L-1$, and whose Fourier transform $\mathcal Fg$ has only isolated zeros.
\end{assumption}

\begin{lemma}
\label{lem:conv_pe}
    Consider ${L\in\N}$, ${L\geq 1}$, and let $g$ satisfy Assumption~\ref{as:1}. A signal $\bar u\in H^{L-1}(\intv,\R^m)$ is persistently exciting of order $L$ if and only if $\tilde u\coloneq \bar u*g$ is persistently exciting of order $L$.
    In this case, if $\col(\bar u,\bar y)\in\B\cap H^{L-1}(\intv, \R^{m+p})$ and $\tilde y\coloneq \bar y*g$, then the Gramians $\Gamma$ in Lemma~\ref{lem:FL} and $\widetilde \Gamma\coloneq \int_\intv \widetilde W(s) \widetilde W^\top(s)\,\mathrm ds$ with $\widetilde W=\col(\tilde u,\dots, \tilde u^{(L-1)}, \tilde y,\dots, \tilde y^{(L-1)})$ have equal images and kernels.
\end{lemma}
The proof of the lemma is given in Appendix~\ref{app:conv_gram}.

The kernels of algebraic differentiators as defined in Appendix~\ref{app:algdiff} satisfy Assumption~\ref{as:1}. The Fourier-domain condition stated in Assumption~\ref{as:1} can be checked using frequency-domain analysis, which also provides insights into the associated noise-attenuation characteristics~\cite{kiltz2013}, \cite{mboup2014}, \cite{othmane2021b}. The condition that $\mathcal{F}g$ has only isolated zeros ensures that the kernel does not annihilate entire frequency bands. Consequently, the convolution with $g$ preserves almost all frequency components of the signal.

Hence, derivatives of $u$ and $y$ estimated using algebraic differentiators can be used to construct $\widetilde{\Gamma}$ as defined in Lemma \ref{lem:conv_pe}. The invariance of image and kernel implies that, in all preceding computations where the Gramian $\Gamma$ appears, it can be replaced by $\widetilde{\Gamma}$ without affecting the results, cf.\ Remark~\ref{rem:replace_gamma}. In particular, the matrix $M$ in \eqref{eq:ddimrep} can be constructed solely from noisy data for a data-based prediction.

\section{Numerical experiments}
\label{sec:experiments}
We first describe how to generate persistently exciting inputs using random splines. Then, we use singular value decomposition (SVD) to compute the matrix $M$ in \eqref{eq:ddimrep}. Finally, we present simulation results.

\subsection{Generation of random persistently exciting inputs}
\label{sec:input_generation}
Let $m$ be the number of components of the input $u$ of system~\eqref{eq:sys}, $L$ the spline degree, $N$ the number of interpolation knots, 
$[t_{\min},t_{\max}]$, $t_{\min}<t_{\max}$,  a time interval, and $A>0$ an amplitude parameter.  
Define a uniform grid 
$t_j = t_{\min} + \tfrac{j}{N-1}(t_{\max}-t_{\min})$ for $j=0,\dots,N-1$, 
and draw a random matrix $Y\in\mathbb{R}^{m\times N}$ with i.i.d.\ entries 
$Y_{i,j}\sim\mathrm{Unif}([-A,A])$.  
For each $i=1,\dots,m$, construct the spline $S_i$ of degree $L$ satisfying 
$S_i(t_j)=Y_{i,j}$.  
The input signal is defined as $\bar u(t)=\col(S_1(t),\dots,S_m(t))$, where $A$ determines its amplitude. 

The next lemma, proved in Appendix~\ref{sec:spline_pe}, justifies this construction.
\begin{lemma}
    \label{lem:spline_pe}
    Suppose $m=1$ and $L<N$. Consider uniform grid points $t_0,\dots,t_{N-1}$ together with i.i.d.\ samples $Y_0,\dots,Y_{N-1}\sim\mathrm{Unif}([-A,A])$, $A>0$. Let $S$ be the spline of degree $L$ satisfying $S(t_i)=Y_i$, $i=0,\dots,N-1$. Then, with probability one, $S$ is persistently exciting of order $L$.
\end{lemma}

\subsection{Regularized pseudo-inverse via SVD}
The matrix $A_3$ in \eqref{eq:ddimrep} may be ill-conditioned. Therefore, we compute its SVD ${A_3 = U \Sigma V^\top}$ and apply Tikhonov regularization with parameter $\lambda_{\text{reg}}$, replacing $A_3^{-1}$ by ${A^{\dagger}_{\text{reg}} = V \operatorname{diag}(s_i/(s_i^2 + \lambda_{\text{reg}})) U^\top}$, where $s_i$ are the singular values of $A_3$ (see \cite{golub2012}).

\subsection{Simulation results and discussion}
Consider
\begin{equation}
		\dot{x}(t) = \begin{bmatrix}
			0 & 1 & 0 \\
			-1 & 0 & 1 \\
			0 & 0 & -2
		\end{bmatrix}\,x(t) + \begin{bmatrix}
		0\\
		0\\
		1
	\end{bmatrix}u(t),\quad
		y(t) = \begin{bmatrix}
			1 & 0 & 1
		\end{bmatrix}x(t),
\end{equation}
with state $x(t) \in \mathbb{R}^3$, input $u(t) \in \mathbb{R}$, and output $y(t) \in \mathbb{R}$.

The persistently exciting input--output trajectory on the interval $[t_{\min}, t_{\max}]$ is sampled uniformly with step size $t_{\mathrm{s}}$, resulting in $N$ samples. A disturbed measurement $\bar y_\eta(t_i)=\bar y(t_i)+\eta_i$, where $\eta_i$, $i\in\N$, is a zero-mean white Gaussian i.i.d. sequence, is used for data-based prediction, and $\bar{y}$ is the output corresponding the persistently exciting input $\bar{u}$. 
The variance of the noise sequence $\eta_i$ is chosen such that the signal-to-noise ratio (SNR) defined as 
\begin{equation}
	\mathrm{SNR} = 10 \log_{10} \left({\sum_{i=1}^N |y(t_i)|^2}\Big /{\sum_{i=1}^N |\eta_i|^2}\right)
	\label{eq:SNR}
\end{equation}
corresponds to a desired value, where $N$ is the number of samples in  $[t_{\min},\,t_{\max}]$.
The signals $\bar u$ and $\bar{y}_\eta$ are then filtered using an algebraic differentiator to approximate the required derivatives. The differentiators are designed and discretized using the toolbox~\cite{othmane2023Tool}, where the trapezoidal rule is employed for numerical integration.  The simulations use the parameter values 
${t_s = 10^{-3}}$, ${\alpha_{\text{diff}}=\beta_{\text{diff}}=8}$, 
${N_{\text{diff}}=0}$, ${T_{\text{diff}} = 84\,t_s}$, 
${\lambda_{\text{reg}} = 10^{-8}}$, 
${t_{\min}=0}$, ${t_{\max}=3\pi}$, 
${N=14}$, $L=4$, and a persistency of excitation order of $L+n=7$, where the subscript \textit{diff} denotes parameters of the differentiator as introduced in Appendix~\ref{app:algdiff}.

 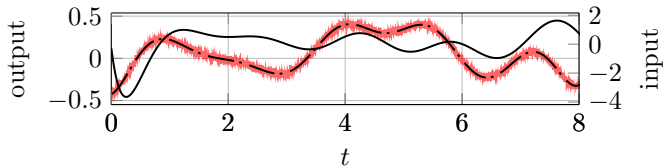
\begin{figure}
	\begin{tikzpicture}
	\pgfplotsset{table/search path={figs/20251122_174404/lambda_1p00em08/noise_2p0000e-02/seed_998994690/}}
	\begin{axis}[
		name=leftaxis,
		width=7.8cm,
		height=2.8cm,
		grid=major,
		xmin=0,
		xmax = 8,
		xlabel={$t$},
		ylabel={output},
		axis x line*=bottom,
		axis y line*=left,
		legend style={at={(0.3,0.38)}, anchor=north, legend columns=2, draw=none},
		]
		\addplot[line width=0.75pt, red, opacity=0.6]
		table[x index=0, y index=1] {pe_output.dat};
		
		\addplot[black, dash pattern={on 10pt off 5pt on 1pt off 5pt}, line width=0.75pt]
		table[x index=0, y index=2] {pe_output.dat};
	\end{axis}
	
	\begin{axis}[
		width=7.8cm,
		height=2.8cm,
		at={(leftaxis.south east)},
		anchor=south east,
		axis y line*=right,
		ylabel={input},
		xmin=0,
		xmax = 8,
		yticklabel style={/pgf/number format/fixed},
		legend style={at={(0.55,0.38)}, anchor=north, legend columns=1, draw=none},
		]
		\addplot[blue, line width=1.5pt]
		table[x index=0, y index=1] {pe_input.dat};
	\end{axis}
\end{tikzpicture}
	\caption{Evolution of input $\bar u$~\legendline{blue, line width=1.5pt}, output 
$\bar y$~\legendline{black, dash pattern={on 10pt off 5pt on 1pt off 5pt}, line width=0.75pt} and disturbed output 
$\bar y_\eta$~\legendline{red, opacity=0.6, line width=0.75pt} with $\mathrm{SNR}= \SI{20.37}{\dB}$.}
	\label{fig:experiment_A}
\end{figure}
 \begin{figure}
%
	\begin{tikzpicture}
		\pgfplotsset{table/search path={figs/20251122_174404/lambda_1p00em08/noise_2p0000e-02/seed_998994690/}}
		\begin{groupplot}[
			group style={
				group size=1 by 2,
				vertical sep=0.45cm,
				xlabels at=edge bottom,
			},
		xlabel=$t$,
			width=8.5cm,
			height=3cm,
			grid=major,
			xmin=0,
			xmax = 16,
			]

			\nextgroupplot[ylabel=signals,	legend style={at={(0.5,1.5)}, anchor=north, legend columns=4, draw=none},]
			\addplot[line width=0.75pt, blue] table[x index=0, y index=1] {traj_state_down.dat};
			\addlegendentry{$\hat{x}_1$}
			\addplot[line width=0.75pt, green] table[x index=0, y index=2] {traj_state_down.dat};
			\addlegendentry{$\hat{x}_2$}
			\addplot[line width=0.75pt, orange] table[x index=0, y index=3] {traj_state_down.dat};
			\addlegendentry{$\hat{x}_3$}

			\addplot[xlabel={$t$},line width=0.75pt, black, dash pattern={on 10pt off 2pt on 1pt off 3pt},] table[x index=0, y index=4] {traj_state_down.dat};
			\addlegendentry{${x}_i$}
			\addplot[line width=0.75pt, black, dash pattern={on 10pt off 2pt on 1pt off 3pt},] table[x index=0, y index=5] {traj_state_down.dat};
			\addplot[line width=0.75pt, black, dash pattern={on 10pt off 2pt on 1pt off 3pt},] table[x index=0, y index=6] {traj_state_down.dat};
			\nextgroupplot[ylabel=input,			legend style={at={(0.7,0.3)}, anchor=north, legend columns=2, draw=none},]
			\addplot[line width=0.75pt, blue] table[x index=0, y index=1] {traj_input_down.dat};
			
		\end{groupplot}
	\end{tikzpicture}
	\caption{Time evolution of the data-based state prediction using the data from Fig. \ref{fig:experiment_A} with respect to the true state and the corresponding input trajectory.}
	\label{fig:results_experiment_A}
 \end{figure}
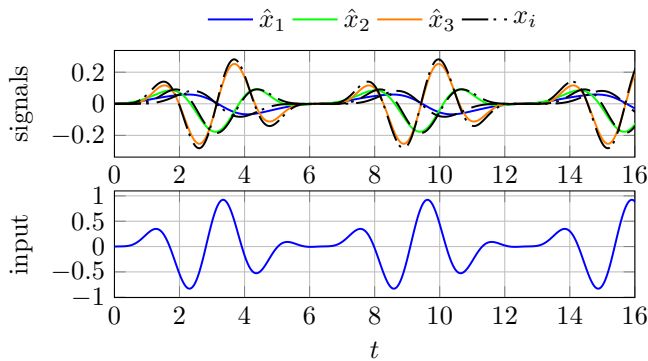

An input–output trajectory is predicted using the data-driven image representation with a chosen $\ell$ defined as a linear combination of trigonometric functions. For benchmarking, an estimate $\hat x$ of the state is reconstructed from the predicted input–output trajectory using the Kalman observability matrix, while a reference state trajectory $x$ is obtained by simulating system~\eqref{eq:sys} with a high-accuracy solver. The initial condition for both trajectories is set to ${x(0)=\hat{x}(0)=\begin{bmatrix}0&0&0\end{bmatrix}{}^\top}$. Let ${X = [x(t_1)\ \cdots\ x(t_N)]}$ and
${\hat X = [\hat x(t_1)\ \cdots\ \hat x(t_N)]}$. For quantitative accuracy assessment, define the relative error in the Frobenius norm
\begin{equation}
  E(x,\hat x) =
  {\bigl\| S^{-1}(\hat X - X) \bigr\|_F}{
        \bigl\| S^{-1} X \bigr\|_F}^{-1},
  \label{eq:relative_error_state}
\end{equation}
where the matrix ${S=\text{diag}(s_1, s_2, \dots, s_n)}$ scales each state $x_i$ using its sample standard deviation.


Fig.~\ref{fig:experiment_A} shows the persistently exciting input and output trajectories for one particular experiment (Experiment I) with $\mathrm{SNR} = \SI{20.37}{\dB}$. The predicted state $\hat{x}$ based on the image representation~\eqref{eq:ddimrep} is shown in Fig.~\ref{fig:results_experiment_A}. The relative error~\eqref{eq:relative_error_state} is $E(x,\hat x) = 0.26$, demonstrating accurate prediction performance, despite a significant noise level.

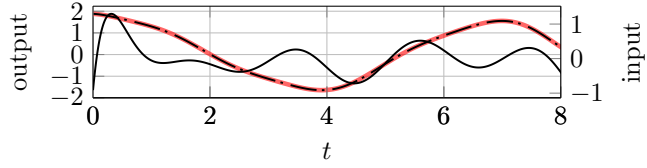
\begin{figure}
	\begin{tikzpicture}
	\pgfplotsset{table/search path={figs/20251122_174404/lambda_1p00em08/noise_1p0000e-04/seed_1871869742/}}
	\begin{axis}[
		name=leftaxis,
		width=7.8cm,
		height=2.8cm,
		grid=major,
		xmin=0,
		xmax = 8,
		xlabel={$t$},
		ylabel={output},
		axis x line*=bottom,
		axis y line*=left,
		legend style={at={(0.183,0.99)}, anchor=north, legend columns=2, fill=white, draw=none},
		]
		\addplot[line width=1.75pt, red, opacity=0.6]
		table[x index=0, y index=1] {pe_output.dat};
		
		\addplot[black, dash pattern={on 10pt off 5pt on 1pt off 5pt}, line width=0.75pt]
		table[x index=0, y index=2] {pe_output.dat};
	\end{axis}
	
	\begin{axis}[
		width=7.8cm,
		height=2.8cm,
		at={(leftaxis.south east)},
		anchor=south east,
		axis y line*=right,
		ylabel={input},
		xmin=0,
		xmax = 8,
		yticklabel style={/pgf/number format/fixed},
		legend style={at={(0.45,0.98)}, anchor=north, legend columns=1, fill=white, draw=none},
		]
		\addplot[blue, line width=1.5pt]
		table[x index=0, y index=1] {pe_input.dat};
	\end{axis}
\end{tikzpicture}
\caption{Evolution of input $\bar u$~\legendline{blue, line width=1.5pt}, output 
$\bar y$~\legendline{black, dash pattern={on 10pt off 5pt on 1pt off 5pt}, line width=0.75pt} and disturbed output 
$\bar y_\eta$~\legendline{red, opacity=0.6, line width=0.75pt}
with $\mathrm{SNR}=\SI{80.44}{\dB}$.
}
	\label{fig:experiment_B}
\end{figure}

\begin{figure}
%

	\begin{tikzpicture}
		\pgfplotsset{table/search path={figs/20251122_174404/lambda_1p00em08/noise_1p0000e-04/seed_1871869742/}}
		\begin{groupplot}[
			group style={
				group size=1 by 1,
				vertical sep=0.45cm,
				xlabels at=edge bottom,
			},
			width=8.5cm,
			height=3cm,
			grid=major,
			xmin=0,
			xmax = 16,
			]

			\nextgroupplot[ylabel=signals, xlabel=$t$,	legend style={at={(0.5,1.5)}, anchor=north, legend columns=4, draw=none},]
			\addplot[line width=0.75pt, blue] table[x index=0, y index=1] {traj_state_down.dat};
			\addlegendentry{$\hat{x}_1$}
			\addplot[line width=0.75pt, green] table[x index=0, y index=2] {traj_state_down.dat};
			\addlegendentry{$\hat{x}_2$}
			\addplot[line width=0.75pt, orange] table[x index=0, y index=3] {traj_state_down.dat};
			\addlegendentry{$\hat{x}_3$}

			\addplot[xlabel={$t$},line width=0.75pt, black, dash pattern={on 10pt off 2pt on 1pt off 3pt},] table[x index=0, y index=4] {traj_state_down.dat};
			\addlegendentry{${x}_i$}
			\addplot[line width=0.75pt, black, dash pattern={on 10pt off 2pt on 1pt off 3pt},] table[x index=0, y index=5] {traj_state_down.dat};
			\addplot[line width=0.75pt, black, dash pattern={on 10pt off 2pt on 1pt off 3pt},] table[x index=0, y index=6] {traj_state_down.dat};
			
		\end{groupplot}
	\end{tikzpicture}
	\caption{Time evolution of the data-based state prediction using the data from  Fig. 
	\ref{fig:experiment_B} with respect to the true state for input trajectory from Fig.\ref{fig:results_experiment_A}.}
    \label{fig:results_experiment_B}
 \end{figure}
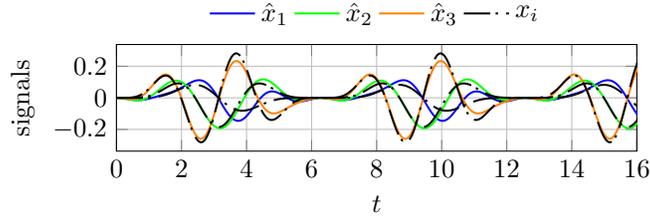

Fig.~\ref{fig:experiment_B} shows the persistently exciting input and output trajectories for another experiment (Experiment II) with a higher $\mathrm{SNR} = \SI{80.44}{\dB}$, yet the predicted state $\hat{x}$ in Fig. \ref{fig:results_experiment_B} yields a larger relative error $E(x,\hat x) = 0.41$. 

\begin{figure}


\begin{tikzpicture}
\pgfplotsset{table/search path={figs/20251122_174404/lambda_1p00em08/}}
\begin{axis}[
    width=7cm,
    height=3.5cm,
    xlabel={SNR in dB},
    ylabel={relative error},
    grid=major,
    grid style={dashed, gray!30},
    ymode=log,
    xmin=13.62,
    xmax=91.85,
    ymin=2.21e-02,
    ymax=1.12e+01,
    legend pos=north east,
    legend cell align=left,
    tick label style={font=\small},
    label style={font=\small},
    title style={font=\bfseries},
    legend style={at={(1.2,1.15)}, anchor=north, legend columns=1, draw=none},
    set layers=axis on top
]
\addlegendimage{empty legend}
\addlegendentry{\hspace{-.0cm}{STD}}
\addplot[
    only marks,
    mark=*,
    mark size=2pt,
    color=cyan,
    fill opacity=0.6,
    draw=gray,
    line width=0.3pt,
] table[x index=0, y index=1] {scatter_data_noise_1p0000em04.dat};
\addlegendentry{$10^{-4}$}

\addplot[
    only marks,
    mark=*,
    mark size=2pt,
    color=black,
    fill opacity=0.6,
    draw=gray,
    line width=0.3pt,
] table[x index=0, y index=1] {scatter_data_noise_5p0000em04.dat};
\addlegendentry{$5\cdot10^{-4}$}

\addplot[
    only marks,
    mark=*,
    mark size=2pt,
    color=orange,
    fill opacity=0.6,
    draw=gray,
    line width=0.3pt,
] table[x index=0, y index=1] {scatter_data_noise_1p0000em03.dat};
\addlegendentry{$10^{-3}$}

\addplot[
    only marks,
    mark=*,
    mark size=2pt,
    color=red,
    fill opacity=0.6,
    draw=gray,
    line width=0.3pt,
] table[x index=0, y index=1] {scatter_data_noise_5p0000em03.dat};
\addlegendentry{$5\cdot10^{-3}$}
\addplot[
    only marks,
    mark=*,
    mark size=2pt,
    color=blue,
    fill opacity=0.6,
    draw=gray,
    line width=0.3pt,
] table[x index=0, y index=1] {scatter_data_noise_1p0000em02.dat};
\addlegendentry{$10^{-2}$}

\addplot[
    only marks,
    mark=*,
    mark size=2pt,
    color=green,
    fill opacity=0.6,
    draw=gray,
    line width=0.3pt,
] table[x index=0, y index=1] {scatter_data_noise_2p0000em02.dat};
\addlegendentry{$2\cdot10^{-2}$}
\pgfplotsextra{
    \begin{pgfonlayer}{axis foreground}
        \coordinate (A) at (axis cs:20.37,0.26);
        \draw[->, thick, black] (axis cs:17,0.7) -- (A);
        \node[anchor=south west, xshift=-0.2cm, yshift=0.1cm] at (axis cs:17,0.7) {\color{black}I};
    \end{pgfonlayer}
}

\pgfplotsextra{
    \begin{pgfonlayer}{axis foreground}
        \coordinate (B) at (axis cs:80.44,0.41);
        \draw[->, thick, black] (axis cs:70,4.5) -- (B);
        \node[anchor=south east, xshift=0.4cm, yshift=-0.1cm] at (axis cs:71,4.0) {\color{black}II};
    \end{pgfonlayer}
}
\end{axis}
\end{tikzpicture}
    
	\caption{Relative error \eqref{eq:relative_error_state} versus SNR \eqref{eq:SNR} for 1500 experiments using  randomly sampled interpolation knots for the input generation and different random noise sequences, colored by noise STD. Letters I and II correspond to results in Figs. \ref{fig:results_experiment_A} and \ref{fig:results_experiment_B}.}
	\label{fig:error_vs_SNR}
\end{figure}
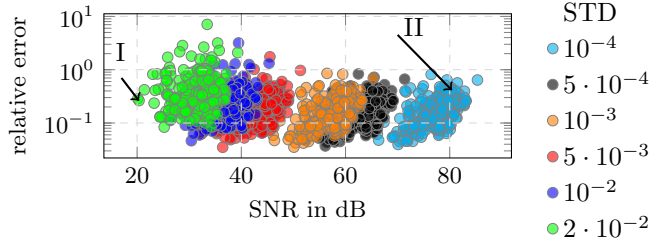

To further analyze performance, we conducted $1500$ experiments with randomly drawn interpolation knots and prescribed noise levels. Fig.~\ref{fig:error_vs_SNR} shows the relative error~\eqref{eq:relative_error_state} versus SNR. While prediction error generally decreases with SNR, indicating robustness, some high-SNR experiments yield larger errors than low-SNR cases, showing that input trajectory quality can be more influential than SNR alone.

\section{Conclusion and outlook}
\label{sec:conclusions}
In this note, we have derived a numerically stable, data-driven method to obtain
an image representation of an unknown continuous-time LTI system, based on a
continuous-time version of Willems et al.'s fundamental lemma. We have employed algebraic differentiators to estimate required derivatives and have proven that the essential excitation conditions are preserved under this derivative approximation scheme. Simulations confirm the effectiveness of the proposed approach, even under significant measurement noise.
Future work should address an extension to systems with multiple outputs and the design of input trajectories that ensure a desired prediction accuracy. Additionally, the obtained results may serve as a basis for direct continuous-time data-driven control methods.

\appendix
\section{Algebraic differentiators}
\label{app:algdiff}

Consider an interval $\mathcal{I}_t=[\timevar-\timewindow,\timevar]\subset \mathbb{R}$, with $\timewindow>0$. By mapping this interval to $[-1,1]$, a square-integrable function $f$ defined on $\mathcal{I}_t$ can be approximated using Jacobi polynomials $\jacopoly{\polyorder}$ of degree $N_{\text{diff}}$, which form an orthogonal basis with respect to the weight 
\begin{equation}
	w(\tau)=\begin{cases}
		(1-\tau)^{\jacoparam}(1+\tau)^{\jacoparambeta},&\tau\in(-1,1)\\
		0,&\text{otherwise}.
	\end{cases}
\end{equation}
The inner product associated with this weight is defined by
\begin{equation}
\langle p , q \rangle_w \;=\; \int_{-1}^{1} p(\tau)\,q(\tau)\,w(\tau)\, \mathrm{d}\tau,
\end{equation}
and the corresponding norm is defined by
$
\|p\|^2_w \;=\; \langle p , p \rangle_w
 \;=\; \int_{-1}^{1} p(\tau)^2\, w(\tau)\, \mathrm{d}\tau.
$

For $\min\{\jacoparam,\jacoparambeta\}>n-1$, $n\geq 0$, a generalized Fourier expansion of order $\polyorder$ yields the estimate 
\begin{equation}\label{eq:estimate_alg_diff}
	\hat f^{(\derivorder)}(\timevar)=\int_{\timevar-\timewindow}^{\timevar} \diffkernel{\derivorder}(\timevar-\tau)f(\tau)\,d\tau=\left(g^{(n)}*f\right)(t)
\end{equation}
for ${f^{(\derivorder)}(t)}$, with
\begin{equation}
		{\diffkernelsymbol(\tau)=\frac{2}{\timewindow}\,w(\nu(\tau))\sum_{j=0}^{\polyorder}\frac{p_j(\vartheta_{\text{diff}})}{\|p_j\|^2_w}p_j(\nu(\tau))},
\end{equation}
 ${\nu(\tau)=1-2\tau/\timewindow}$,  window length $\timewindow$,  Jacobi polynomial $p_j=\jacopoly{j}$ with parameters $\jacoparam$ and $\jacoparambeta$,  expansion order $\polyorder$, and the scalar $\vartheta_\text{diff}$ parameterizing the approximation delay (see, e.g., \cite{othmane2021b} for more details). Notice that $g^{(n)}*f=g*f^{(n)}$.

\section{Convoluted Gramians}
Before we prove Lemma~\ref{lem:conv_pe} in this section, we state and prove the next auxiliary result.
\label{app:conv_gram}
\begin{lemma}
\label{lem:conv_gram}
    Let $g$ satisfy Assumption~\ref{as:1}. Consider $w\in H^{L-1}(\intv, \R^d)$. Then for $W:=\col(w,\dots,w^{(L-1)})$ and $W_g:=\col(w*g,\dots, (w*g)^{(L-1)})$ the Gramians
    \begin{equation}
        \Gamma \coloneq \int_\intv W(s)W(s)^\top\,\mathrm ds,\quad \Gamma_g \coloneq \int_\intv W_g(s)W_g(s)^\top\,\mathrm ds
    \end{equation}
    have the same image and kernel. In particular, $\eta\in\ker\Gamma = \ker\Gamma_g$ if and only if $W^\top \eta = W_g^\top \eta = 0$ a.e.
\end{lemma}
\begin{proof}
    Observe that $W_g = W*g = g*W$. We show $\ker \Gamma = \ker\Gamma_g$. Note that $\Gamma \eta=0$ if and only if $\eta^\top \Gamma \eta=0$, which is equivalent to $W^\top \eta=0$ a.e.; similarly for $\Gamma_g$ and $W_g$. Let $\eta\in\ker \Gamma$ and, therefore, $W_g^\top \eta = g*(W^\top\eta) =0$ a.e. This implies $\eta\in\ker\Gamma_g$. Let $\eta\in\ker\Gamma_g$ and, therefore, $W_g^\top \eta = 0$ a.e. Applying the Fourier transform one sees $$0=\mathcal F (W_g^\top \eta) = \mathcal F(g*W^\top\eta) =\mathcal F g \cdot \mathcal F (W^\top \eta)\quad \text{a.e.}$$
    Since $\mathcal F g$ has only isolated zeros this implies $\mathcal F (W^\top \eta)=0$ a.e.\ and, hence, $W^\top \eta=0$ a.e. This yields $\eta\in\ker \Gamma$.

    Since both Gramians are symmetric matrices, we find 
    \begin{equation*}
        \operatorname{im}\Gamma = (\ker\Gamma)^\bot = (\ker\Gamma_g)^\bot = \operatorname{im}\Gamma_g.\qedhere
    \end{equation*}
\end{proof}

Now we turn to Lemma~\ref{lem:conv_pe}. We apply Lemma~\ref{lem:conv_gram} with $w=\bar u$. By definition, $\bar u$ is persistently exciting of order $L$ if and only if $W^\top \eta = 0$ implies $\eta=0$, which by Lemma~\ref{lem:conv_gram} is equivalent to $\ker \Gamma = \{0\}$; a similar statement is true for $\tilde u=\bar u*g=w_g$ and $\Gamma_g$. This shows the first assertion in Lemma~\ref{lem:conv_pe}. The second claim in Lemma~\ref{lem:conv_pe} follows immediately by applying Lemma~\ref{lem:conv_gram} with $w=\col(\bar u,\bar y)$. 

\section{Persistently exciting splines}
\label{sec:spline_pe}
In this section, we prove Lemma~\ref{lem:spline_pe}, by showing that the complementary event $\Omega_0$, that the spline $S$ is not persistently exciting of order $L$, has probability zero. Assume that $S$ is not persistently exciting of order $L$. It is not difficult to see that this implies that the restriction of $S$ to any subinterval is not persistently exciting of order $L$ either.

On each interval $\intv_i:=(t_i,t_{i+1})$, the spline $S$ coincides with a polynomial $p_i\in\R[s]$ of degree at most $L$. Since polynomials of degree $L$ are persistently exciting of order $L$ (see \cite[Ex.~18]{Schmitz2024a}), in fact, $\deg p_i \leq L-1$ for $i=0,\dots, N-1$. Note that the spline $S$ is $(L-1)$-times continuously differentiable. In particular,
$p_i^{(j)}(t_i) = S^{(j)}(t_i) = p_{i+1}^{(j)}(t_i)$
for $i=0,\dots,N-2$ and $j=0,\dots, L-1$. Because a degree-$(L-1)$ polynomial is uniquely determined by its first $L$ derivatives at a single point, all polynomials $p_i$ coincide. Therefore, $S$ is a polynomial on $[t_0,t_{N-1}]$ with $\deg S\leq L-1$. This implies $\Omega_0\subset \Omega\cap [-A,A]^N$, where
\begin{equation}\Omega=
    \left\{\left[\begin{smallmatrix}
Y_0 \\ \vdots \\ Y_{N-1}
\end{smallmatrix}\right]\in \R^{N}\,\middle|\,\begin{aligned} &\exists\, p\in\R[s],\ \deg p \leq L-1,\\&p(t_i)=Y_i,\ i=0,\dots, N-1 \end{aligned}\right\}.
\end{equation}
We show that $\Omega$ has Lebesgue measure zero. Note that the space of real polynomials of degree at most $L-1$ is isomorphic to $\R^{L}$. Let $\Pi$ denote the corresponding isomorphism and $\delta$ be the point-evaluation operator $\delta:\R[s]\rightarrow \R^{N}$, $p\mapsto \begin{bmatrix}
p(t_0) & \dots & p(t_{N-1})
\end{bmatrix}{}^\top$. The dimension of $\Omega=\delta(\Pi(\R^{L}))$ is bounded above by $L$, which by $L<N$ makes $\Omega$ a proper linear subspace of $\R^N$ and, thus, it is a zero-set with respect to the $N$-dimensional Lebesgue measure $\lambda\hspace{-1.25ex}\lambda$. Consequently, the probability of the event $\Omega_0$ is bounded by
\begin{equation*}
\mathbf P(\Omega_0) \leq \int_{\Omega\cap [-A,A]^N} (2A)^{-1}\,\mathrm d\lambda\hspace{-1.25ex}\lambda = 0.
\end{equation*} 
\bibliographystyle{siam}
\bibliography{ref}

@article{dorfler2022bridging,
  title={Bridging direct and indirect data-driven control formulations via regularizations and relaxations},
  author={D{\"o}rfler, Florian and Coulson, Jeremy and Markovsky, Ivan},
  journal={IEEE Trans. Autom. Control},
  volume={68},
  pages={883--897},
  year={2022},
  publisher={IEEE}
}

@Article{mboup2014,
	author    = {M. Mboup and S. Riachy},
	journal   = {{IFAC} {P}roc. Volumes},
	title     = {A Frequency Domain Interpretation of the Algebraic Differentiators},
	year      = {2014},
	pages     = {9147--9151},
	volume    = {47},
	publisher = {Elsevier BV},
}

@Book{golub2012,
	author    = {Golub, G. H. and Van Loan, C. F.},
	publisher = {Johns Hopkins University Press},
	title     = {Matrix Computations},
	year      = {2013},
	address   = {Baltimore},
}

@Article{mboup2009a,
	author    = {M. Mboup and C. Join and M. Fliess},
	journal   = {Numer. Algorithms},
	title     = {Numerical differentiation with annihilators in noisy environment},
	year      = {2009},
	pages     = {439--467},
	volume    = {50},
	publisher = {Springer Nature},
}

@InProceedings{kiltz2013,
	author    = {L. {Kiltz} and J. {Rudolph}},
	booktitle = {52nd IEEE Conf. Decis. Control (CDC)}, 
	title     = {Parametrization of algebraic numerical differentiators to achieve desired filter characteristics},
	year      = {2013},
	pages     = {7010--7015},
}

@Article{othmane2023Tool,
	author  = {Othmane, A. and Rudolph, J.},
	journal = {at - Autom.},
	title   = {AlgDiff: an open source toolbox for the design, analysis and discretisation of algebraic differentiators},
	year    = {2023},
	pages   = {612--623},
	volume  = {71},
	doi     = {doi:10.1515/auto-2023-0035},
}

@Article{othmane2021b,
	author  = {Othmane, A. and Kiltz, L. and Rudolph, J.},
	journal = {Int. J. Syst. Sci.},
	title   = {Survey on algebraic numerical differentiation: Historical developments, parametrization, examples, and applications},
	year    = {2022},
	pages   = {1848--1887},
	volume  = {53},
	doi     = {10.1080/00207721.2022.2025948},
}

@Article{Schmitz2024a,
  author    = {Schmitz, Philipp and Faulwasser, Timm and Rapisarda, Paolo and Worthmann, Karl},
  journal   = {Systems \& Control Letters},
  title     = {A continuous-time fundamental lemma and its application in data-driven optimal control},
  year      = {2024},
  volume    = {194, 105950},
  doi       = {10.1016/j.sysconle.2024.105950},
  publisher = {Elsevier BV},
}

@Article{Beelen1988,
  author    = {Beelen, T. and van Dooren, P.},
  journal   = {SIAM J. Matrix Anal. Appl.},
  title     = {A Pencil Approach for Embedding a Polynomial Matrix into a Unimodular matrix},
  year      = {1988},
  pages     = {77--89},
  volume    = {9},
  doi       = {10.1137/0609006},
  publisher = {Society for Industrial & Applied Mathematics (SIAM)},
}

@Article{oarua1997improved,
  author    = {Oar{\u{a}}, C and Van Dooren, Paul},
  journal   = {Systems \& {C}ontrol {L}etters},
  title     = {An improved algorithm for the computation of structural invariants of a system pencil and related geometric aspects},
  year      = {1997},
  pages     = {39--48},
  volume    = {30},
  publisher = {Elsevier},
}

@Article{Willems2005,
  author     = {Willems, Jan C. and Rapisarda, Paolo and Markovsky, Ivan and De Moor, Bart L. M.},
  journal    = {Systems \& Control Letters},
  title      = {A note on persistency of excitation},
  year       = {2005},
  pages      = {325--329},
  volume     = {54},
  doi        = {10.1016/j.sysconle.2004.09.003},
  publisher  = {Elsevier}
}

@Book{Willems1997,
  author    = {Willems, Jan C. and Polderman, Jan W.},
  publisher = {Springer},
  title     = {Introduction to mathematical systems theory: {A} behavioral approach},
  year      = {1997}
}

@InProceedings{Coulson2019,
  author       = {Coulson, Jeremy and Lygeros, John and D{\"o}rfler, Florian},
  booktitle    = {18th European Control Conference (ECC)},
  title        = {Data-enabled predictive control: In the shallows of the {D}ee{PC}},
  year         = {2019},
  pages        = {307--312},
  doi          = {10.23919/ecc.2019.8795639},
  journal      = {2019 18th European Control Conference (ECC)}
}

@Article{Persis2020,
  author    = {De Persis, Claudio and Pietro Tesi},
  journal   = {IEEE Trans. Autom. Control},
  title     = {Formulas for Data-Driven Control: Stabilization, Optimality, and Robustness},
  year      = {2020},
  pages     = {909--924},
  volume    = {65},
  doi       = {10.1109/tac.2019.2959924},
  publisher = {Institute of Electrical and Electronics Engineers (IEEE)},
}

@Article{Markovsky2021,
  author    = {Ivan Markovsky and Florian Dörfler},
  journal   = {Annual Reviews in Control},
  title     = {Behavioral systems theory in data-driven analysis, signal processing, and control},
  year      = {2021},
  pages     = {42--64},
  volume    = {52},
  doi       = {10.1016/j.arcontrol.2021.09.005},
  publisher = {Elsevier BV},
}

@Article{Faulwasser2023,
  author     = {Faulwasser, Timm and Ou, Ruchuan and Pan, Guanru and Schmitz, Philipp and Worthmann, Karl},
  journal    = {Annual Reviews in Control},
  title      = {Behavioral theory for stochastic systems? {A} data-driven journey from {W}illems to {W}iener and back again},
  year       = {2023},
  pages      = {92--117},
  doi        = {10.1016/j.arcontrol.2023.03.005},
  publisher  = {Elsevier},
}

@InProceedings{Lopez_2022,
  author    = {Lopez, Victor G. and M\"uller, Matthias A.},
  booktitle = {61st IEEE Conf. Decis. Control (CDC)},
  title     = {On a Continuous-Time Version of {W}illems’ Lemma},
  year      = {2022},
  pages     = {2759--2764},
  doi       = {10.1109/cdc51059.2022.9992347},
  journal   = {2022 IEEE 61st Conference on Decision and Control (CDC)},
}

@Article{Lopez2024,
  author    = {Victor G. Lopez and Matthias A. Müller and Paolo Rapisarda},
  journal   = {IEEE Control Systems Letters},
  title     = {An Input-Output Continuous-Time Version of {W}illems’ Lemma},
  year      = {2024},
  pages     = {916--921},
  volume    = {8},
  doi       = {10.1109/lcsys.2024.3406057},
  publisher = {Institute of Electrical and Electronics Engineers (IEEE)},
}

@article{Rapisarda2023,
title = {A “fundamental lemma” for continuous-time systems, with applications to data-driven simulation},
journal = {{S}ystems \& {C}ontrol {L}etters},
volume = {179, 105603},
year = {2023},
doi = {https://doi.org/10.1016/j.sysconle.2023.105603},
author = {P. Rapisarda and M.K. Çamlibel and H.J. {van Waarde}},
}

@article{Schmitz2026,
    author = {P. Schmitz and K. Worthmann and T. Faulwasser and P. Rapisarda},
    title = {Data-driven continuous-time optimal control:
A unified framework using orthogonal functions},
    journal = {to appear in Math.\ Control Signals Syst.},
}
\end{document}